\documentclass{article}%
\usepackage{amsmath}
\usepackage{amsfonts}
\usepackage{amssymb}
\usepackage{graphicx}%
\providecommand{\U}[1]{\protect\rule{.1in}{.1in}}
\newtheorem{theorem}{Theorem} [section]

\newtheorem{conjecture}[theorem]{Conjecture}

\newtheorem{proposition}[theorem]{Proposition}

\newenvironment{proof}[1][Proof]{\noindent\textbf{#1.} }{\ \rule{0.5em}{0.5em}}
\begin{document}

\title{A Simple Proof of an Inequality Connecting the Alternating Number of
Independent Sets and the Decycling Number}
\author{Vadim E. Levit\\Ariel University Center of Samaria, Ariel, Israel\\levitv@ariel.ac.il
\and Eugen Mandrescu\\Holon Institute of Technology, Holon, Israel\\eugen\_m@hit.ac.il}
\date{}
\maketitle

\begin{abstract}
If $s_{k}$ denotes the number of independent sets of cardinality $k$ and
$\alpha(G)$ is the size of a maximum independent set in graph $G$, then
$I(G;x)=s_{0}+s_{1}x+...+s_{\alpha(G)}x^{\alpha(G)}$ is the
\textit{independence polynomial} of $G$ \cite{GutHar}.

In this paper we provide an elementary proof of\ the inequality\textit{
}$\left\vert I(G;-1)\right\vert \leq2^{\varphi(G)}$, where $\varphi(G)$ is the
decycling number of $G$.


\end{abstract}

\section{Introduction}

Throughout this paper $G=(V,E)$ is a finite, undirected, loopless and without
multiple edges graph, with vertex set $V=V(G)$ and edge set $E=E(G)$. By $G-W$
we mean the subgraph induced by $V-W$. The set $N(v)=\{w:w\in V$
\ \textit{and} $vw\in E\}$ \textit{neighborhood} of the vertex $v\in V$, and
$N[v]=N(v)\cup\{v\}$. A \textit{leaf} is a vertex having a unique neighbor.

A set of pairwise non-adjacent vertices is called \textit{independent}. The
\textit{independence number }of $G$, denoted by $\alpha(G)$, is the
cardinality of a maximum independent set.

If $G$ has $s_{k}$ independent sets of size $k$, then
\[
I(G;x)=s_{0}+s_{1}x+s_{2}x^{2}+...+s_{\alpha(G)}x^{\alpha(G)}%
\]
is known as the \textit{independence polynomial} of $G$ \cite{GutHar}. Some
properties of the independence polynomial are presented in
\cite{AlMalSchErdos,Brown,ChudSeymour,LeMan04,LeMa04b,LeMa04e}.

The value of a graph polynomial at a specific point can give sometimes a very
surprising information about the structure of the graph \cite{BalBolCutPeb}.
In the case of independence polynomials, let us notice that:

\begin{itemize}
\item $I(G;1)=s_{0}+s_{1}+s_{2}+...+s_{\alpha}$ equals the number of
independent sets of $G$. It is known as the \textit{Fibonacci number} of $G$
\cite{KTiWagZieg,PedVest,ProdTichy}.

\item $I(G;-1)=s_{0}-s_{1}+s_{2}-...+(-1)^{\alpha}s_{\alpha}$ is equal to
difference of the numbers of independent sets of even and odd sizes. It is
known as the \textit{alternating number of independent sets} \cite{BSN07}. The
value of $\left\vert I(G;-1)\right\vert $ can be any non-negative integer. For
instance, $\left\vert I(K_{\alpha,\alpha,...,\alpha};-1)\right\vert =n-1$,
where $K_{\alpha,\alpha,...,\alpha}$ is the complete $n$-partite graph.

\item $I(G;-1)=-\widetilde{\chi}\left(  Ind(G)\right)  $, where $\widetilde
{\chi}\left(  \Sigma\right)  $ is the \textit{reduced Euler characteristic }of
the abstract simplicial complex $\Sigma$. Recall that an abstract simplicial
complex on a finite vertex set $\Sigma_{0}$ is a subset $\Sigma$ of
$2^{\Sigma_{0}}$ satisfying: $\{v\}\in\Sigma$ for every $v\in\Sigma_{0}$, and
$A\subseteq B\in\Sigma$ implies $A\in\Sigma$. The elements of $\Sigma$ are
\textit{faces} and the dimension of a face $A$ is $\left\vert A\right\vert
-1$. For a simplicial complex with $s_{i}$ faces of dimension $i-1$, the
\textit{reduced Euler characteristic} equals $-s_{0}+s_{1}-s_{2}+s_{3}-...$.
The family $Ind(G)$ of all independent sets of a graph $G=(V,E)$ forms a
simplicial complex on $V$, called the \textit{independence complex} of $G$
\cite{Jonsson08}.
\end{itemize}

The \textit{cyclomatic number} $\nu(G)$ of the graph $G$ is the dimension of
the \textit{cycle space} of $G$, i.e., the dimension of the linear space
spanned by the edge sets of all the cycles of $G$. The \textit{decycling
number }\cite{Beineke} (or the \textit{feedback vertex number} \cite{Ueno})
$\varphi(G)$ of a graph $G$ is the minimum number of vertices that need to be
removed in order to eliminate all its cycles. While $\nu(G)$ can be easily
computed, since $\nu(G)=\left\vert E(G)\right\vert -\left\vert V(G)\right\vert
+p$, where $p$ is the number of connected components of $G$, it is known that
to compute $\varphi(G)$ is an \textbf{NP}-complete problem \cite{Karp}. It is
clear that $\varphi(G)\leq\nu(G)$ holds for every graph $G$.

The inequality $\left\vert I(G;-1)\right\vert \leq2^{\nu(G)}$ has been
established in \cite{LevMan09}, while a stronger result, namely, $\left\vert
I(G;-1)\right\vert \leq2^{\varphi(G)}$ has been proved in \cite{Engstrom}.

In this paper we provide a simple proof of the inequality $\left\vert
I(G;-1)\right\vert \leq2^{\varphi(G)}$ using only elementary arguments.

\section{Results}

\begin{proposition}
\cite{GutHar}\label{prop1} If $v\in V(G)$, then $I(G;x)=I(G-v;x)+x\cdot
I(G-N[v];x)$.
\end{proposition}

\begin{theorem}
For any graph $G$ the alternating number of independent sets is bounded as
follows
\[
\left\vert I(G;-1)\right\vert \leq2^{\varphi(G)},
\]
where $\varphi(G)$ is the decycling number of $G$.
\end{theorem}

\begin{proof}
We establish the inequality by induction on $\varphi(G)$.

\begin{itemize}
\item If $\varphi(G)=0$, then $G$ is a forest, and we have to show that
$\left\vert I(G;-1)\right\vert \leq1$.

We proceed by mathematical induction on $n=\left\vert V(G)\right\vert $.

For $n=0,I(G;x)=1$ and $I(G;-1)=1$, while for $n=1,I(G;x)=1+x$ and
$I(G;-1)=0$. Suppose that $G$ is a forest with $\left\vert V(G)\right\vert
=n\geq2.$

If $G$ has no leaves, then $I(G;x)=(1+x)^{n}$ and $I(G;-1)=0$. Otherwise, let
$v$ be a leaf of $G$ and $N(v)=\{u\}$. According to Proposition \ref{prop1} we
obtain that
\[
I(G;x)=I(G-u;x)+x\cdot I(G-N[u];x)=(1+x)\cdot I(G-\{u,v\};x)+x\cdot
I(G-N[u];x).
\]
Hence, by the induction hypothesis, we finally get
\[
\left\vert I(T;-1)\right\vert =\left\vert (-1)\cdot I(T-N[u];-1)\right\vert
\leq1.
\]

\item Assume that the statement is true for graphs with the decycling number
$\varphi(G)\leq k$.

Let $G$ be a graph with $\varphi(G)=k+1$. Clearly, there exists some $v\in
V(G)$, such that $\varphi(G-v)<\varphi(G)$. According to Proposition
\ref{prop1}, we get:
\[
I(G;-1)=I(G-v;-1)-I(G-N[v];-1).
\]
By the induction hypothesis, it assures that
\[
\left\vert I(G;-1)\right\vert \leq\left\vert I(G-v;-1)\right\vert +\left\vert
I(G-N[v];-1)\right\vert \leq2\cdot2^{k}=2^{\varphi(G)},
\]
and this completes the proof.
\end{itemize}
\end{proof}

Notice that if $G=qK_{3}$, then $I(G;x)=(1+3x)^{q}$ and hence,
$I(G;-1)=(-2)^{\varphi(G)}$.

\begin{conjecture}
For every positive integer $k$ and each integer $q$ such that $\left\vert
q\right\vert \leq2^{k}$, there is a graph $G$ with $\varphi(G)=k$ and
$I(G;-1)=q$.
\end{conjecture}

\section{Acknowledgements}

The authors thank Alexander Engstr\"{o}m for drawing their attention to the
fact that some of their findings from \cite{LevMan09} may be translated to
known results in combinatorial topology.


\begin{thebibliography}{99}                                                                                               %
\bibitem {AlMalSchErdos}Y. Alavi, P. J. Malde, A. J. Schwenk, P. Erd\"{o}s,
\emph{The vertex independence sequence of a graph is not constrained},
Congressus Numerantium 58 (1987) 15-23.

\bibitem {BalBolCutPeb}P. N. Balister, B. Bollob\'{a}s, J. Cutler, L. Pebody,
\emph{The interlace polynomial of graphs at }$-1$, European Journal of
Combinatorics 23 (2002) 761-767.

\bibitem {Beineke}L. W. Beineke, R. C. Vandell, \emph{Decycling graphs},
Journal of Graph Theory \textbf{25} (1997) 59-77.

\bibitem {Brown}J. I. Brown, K. Dilcher, R. J. Nowakowski, \emph{Roots of
independence polynomials of well-covered graphs}, Journal of Algebraic
Combinatorics \textbf{11} (2000) 197-210.

\bibitem {BSN07}M. Bousquet-{M}\'{e}lou, S. Linusson, E. Nevo, \emph{On the
independence complex of square grids}, Journal of Algebraic Combinatorics
\textbf{27} (2008) 423-450.

\bibitem {ChudSeymour}M. Chudnovsky, P. Seymour, \emph{The roots of the
independence polynomial of a clawfree graph}, Journal of Combinatorial Theory
B \textbf{97} (2007) 350--357.

\bibitem {Engstrom}A. Engstr\"{o}m, \emph{Upper bounds on the Witten index for
supersymmetric laticce models by discrete Morse theory}, European Journal of
Combinatorics \textbf{30} (2009) 429-438.

\bibitem {GutHar}I. Gutman, F. Harary, \emph{Generalizations of the matching
polynomial}, Utilitas Mathematica \textbf{24} (1983) 97-106.

\bibitem {Jonsson08}J. Jonsson, \emph{Simplicial Complexes of Graphs}, Lecture
Notes in Mathematics \textbf{1928}, Springer, 2008, 378 pp.

\bibitem {Karp}R. Karp, \emph{Reducibility among combinatorial problems}, in
\textit{Complexity of Computer Computations} (eds. R. E. Miller and J. W.
Thatcher), Plenum Press 1972, 85-103.

\bibitem {KTiWagZieg}A. Knopfmachera, R. F. Tichy, S. Wagner, V. Ziegler,
\emph{Graphs, partitions and Fibonacci numbers}, Discrete Applied Mathematics
\textbf{155} (2007) 1175-1187.

\bibitem {LeMan04}V. E. Levit, E. Mandrescu, \emph{The independence polynomial
of a graph - a survey}, Proceedings of the $1^{st}$ International Conference
on Algebraic Informatics, Aristotle University of Thessaloniki, Greece, 20-23
October, 2005, pp. 233-254.

\bibitem {LeMa04b}V. E. Levit, E. Mandrescu, \emph{Independence polynomials of
well-covered graphs: Generic counterexamples for the unimodality conjecture},
European Journal of Combinatroics \textbf{27} (2006) 931-939.

\bibitem {LeMa04e}V. E. Levit, E. Mandrescu, \emph{On the roots of
independence polynomials of almost all very well-covered graphs}, Discrete
Applied Mathematics \textbf{156} (2008) 478-491.

\bibitem {LevMan09}V. E. Levit, E. Mandrescu, \emph{The independence
polynomial of a graph at -1}, E-print, arXiv:0904.4819v1 [math.CO]

\bibitem {PedVest}A. S. Pedersen, P. D. Vestergaard, \emph{The number of
independent sets in unicyclic graphs}, Discrete Applied Mathematics
\textbf{152} (2005) 246-256.

\bibitem {ProdTichy}H. Prodinger, R. F. Tichy, \emph{Fibonacci numbers of
graphs}, Fibonacci Quart. \textbf{20} (1982) 16-21.

\bibitem {Ueno}S. Ueno, Y. Kajitani, S. Gotoh, \emph{On the nonseparating
independent set problem and feedback set problem for graphs with no vertex
degree exceeding three}, Discrete Mathematics \textbf{72} (1988) 355-360.
\end{thebibliography}
\end{document}